# A novel methodology of weighted residual for nonlinear computations

W. Chen

Permanent mail address: Dr. Wen CHEN, P. O. Box 2-19-201, Jiangshu University of Science & Technology, Zhenjiang City, Jiangsu Province 212013, P. R. China

Present mail address (as a JSPS Postdoctoral Research Fellow): Apt.4, West 1$^{st}$ floor, Himawari-so, 316-2, Wakasato-kitaichi, Nagano-city, Nagano-ken, 380-0926, JAPAN

E-mail: chenw@homer.shinshu-u.ac.jp

Permanent email: chenwwhy@hotmail.com

**Abstract**

One of strengths in the finite element (FE) and Galerkin methods is their capability to apply weak formulations via integration by parts, which leads to elements matching at lower degree of continuity and relaxes requirements of choosing basis functions. However, when applied to nonlinear problems, the methods of this type require a great amount of computing effort of repeated numerical integration. It is well known that the method of weighted residual is the very basis of various popular numerical techniques including the FE and Galerkin methods. This paper presents a novel methodology of weighted residual for nonlinear computation with objectives to avoid the above-mentioned shortcomings. It is shown that the presented nonlinear formulations of the FE and Galerkin methods can be expressed in the Hadamard product form as in the collocation and finite difference methods. Therefore, the recently developed SJT product approach can be applied in the evaluation of the Jacobian matrix of these nonlinear formulations. This also provides possibility to introduce the nonlinear uncoupling technique to the FE and Galerkin nonlinear computing. Furthermore, the present scheme of weighted residuals also greatly eases the use of the least square and boundary element methods to nonlinear problems.

1. **Introduction**

Today numerical computations have been playing a key role in nearly all science and engineering area. It is recognized that the numerical techniques for solving linear problems have been well established. On the other hand, although the rapid advance in the design of computer has been achieved in recent years, it is still a tough challenge to numerically simulate the large complicated nonlinear systems such as weather prediction and turbulent fluid. An inordinate amount of computing time and storage raised due to nonlinearity often prohibits such calculations. Therefore, scientists and engineers have more expectations in algorithms than in computer. The development of novel computational schemes is therefore in much demand to increase the scale, accuracy and reliability of nonlinear computations.

The finite element and Galerkin methods are currently the standard numerical technique in use to solve various nonlinear problems. The methods retain the advantages of weak formulations, which lowers the continuity requirements of matching elements and permits to use simple basis functions. However, these methods demand a great amount of numerical integration effort in updating Jacobian matrix of each Newton-Raphson iteration. It was reported that numerical integration often occupied a large proportion of CPU time in the FE and Galerkin solution of large nonlinear systems [1]. On the other hand, the least square and boundary element (BE) methods encounter large obstacle to be applied to general nonlinear problems.

Recently, the present author [2, 3] applied the Hadamard product to express the nonlinear formulations of the finite difference (FE) and collocation (pseudo-spectral) methods in an explicit matrix form. Moreover, the SJT product was therein introduced to evaluate the Jacobian matrix efficiently and accurately. A nonlinear uncoupling technique was also developed by means of the Hadamard and SJT product approach [3, 4]. The simplicity and efficiency of such technique were well demonstrated through testing some benchmark problems. In contrast, when the FE, Galerkin, BE, and least square method are applied to nonlinear problems, the corresponding nonlinear analogue formulations can not be expressed as simple and explicit matrix form of the Hadamard product [3]. Therefore, the SJT product and uncoupling techniques developed in [2-4] can not be extended to the nonlinear computations of these numerical schemes.

It is well known that all these numerical techniques have their roots on the weighted residuals method. The objective of this paper is to introduce a novel methodology of weighted residual, which holds the merits of these numerical techniques but overcomes the above-mentioned weaknesses. The remaining of this paper is divided into three sections. Section 2 illustrates an original weighted residuals formulation for nonlinear analogue. It is shown that the presented nonlinear formulations can be expressed in the Hadamard product form. In section 3, the recently developed SJT approach is used to efficiently calculate the accurate Jacobian matrix of these nonlinear formulations, and a approach of uncoupling coupled nonlinear modeling is presented. Finally, some remarks are given in section 4. The present work can be regarded as a continuation of research effort in nonlinear computations [2-4].

## 2. A novel methodology of weighted residuals

The method of weighted residuals (MWR) can be recognized the origin of almost all popular numerical techniques [5, 6]. Consider the differential equations of the form

$$\psi\{u\} - f = 0, \text{ in } \Omega \qquad (1)$$

with the following boundary conditions

$$u = \bar{u}, \text{ on } \Gamma_1 \qquad (2a)$$

$$q = \partial u/\partial n = \bar{q}, \quad \text{on} \quad \Gamma_2, \tag{2b}$$

where n is the outward normal to the boundary, $\Gamma=\Gamma_1+\Gamma_2$, and the upper bars indicate known values. More complex boundary conditions can be easily included but they will not be considered here for the sake of brevity. In the MWR, the desired function u in the differential governing equation is first approximated by a set of linearly independent basis functions $\phi_k(x)$, such that

$$u = \hat{u} = \sum_{k=1}^{n} c_k \phi_k, \tag{3}$$

where $c_k$'s are the unknown parameters. In the Galerkin and FE methods, the basis functions are usually chosen so as to satisfy certain given conditions such as the boundary conditions and the degree of continuity. In addition, these basis functions should be complete.

Substituting equation (3) into equation (1) produces an error, which is called the residual, namely,

$$\psi\{\hat{u}\} - f = R \neq 0. \tag{4}$$

This error or residual R is forced to be zero in the average sense by setting weighted integral of the residuals equal to zero, namely,

$$\int_\Omega [\psi\{\hat{u}\} - f] W_j d\Omega = \int_\Omega R W_j d\Omega = 0, \quad j=1,2,...,N, \tag{5}$$

where $W_j$'s are weighting functions. The differences among weighting functions and basis functions give rise to different numerical techniques such the Galerkin, least square, finite element, boundary element, moments, spectral methods, finite difference and collocation methods.

This paper places its emphasis on the nonlinear computations. Let us consider the quadratic nonlinear operator of the form:

$$p(u)r(u) + L(u) = f, \qquad (6)$$

where p(u), r(u) and L(u) are linear differential operators, f is the constant. The traditional scheme of weighted residuals approximates equation (6) by

$$\int_\Omega [p(\hat{u})r(\hat{u}) + L(\hat{u}) - f]W_j d\Omega = 0, \quad j=1,2,\ldots,N. \qquad (7)$$

Here we present an innovative scheme

$$\int_\Omega p(\hat{u})W_j d\Omega \int_\Omega r(\hat{u})W_j d\Omega + \int_\Omega [L(\hat{u}) - f]W_j d\Omega = 0, \quad j=1,2,\ldots,N \qquad (8)$$

different from the standard equation (7). It is noted that the key distinction of equations (7) and (8) lies in that latter weights linear operators p( ) and r( ) separately. The nonlinear operators are here understood certain combinations of linear operators. Before further development, we first introduce the concepts of the Hadamard matrix product, power and function.

**Definition 2.1** Let matrices $A=[a_{ij}]$ and $B=[b_{ij}] \in C^{N \times M}$, the Hadamard product of matrices is defined as $A \circ B = [a_{ij} b_{ij}] \in C^{N \times M}$. where $C^{N \times M}$ denotes the set of N×M real matrices.

**Definition 2.2** If matrix $A=[a_{ij}] \in C^{N \times M}$, then $A^{\circ q}=[a_{ij}^q] \in C^{N \times M}$ is defined as the Hadamard power of matrix A, where q is a real number. Especially, if $a_{ij} \neq 0$, $A^{\circ(-1)}=[1/a_{ij}] \in C^{N \times M}$ is defined as the Hadamard inverse of matrix A. $A^{\circ 0}=11$ is defined as the Hadamard unit matrix in which all elements are equal to unity.

**Definition 2.3** If matrix $A=[a_{ij}]\in C^{N\times M}$, then the Hadamard matrix function $f^{\circ}(A)$ is defined as $f^{\circ}(A)=[f(a_{ij})]\in C^{N\times M}$.

**Theorem 2.1**: letting A, B and $Q\in C^{N\times M}$, then

1. $A\circ B=B\circ A$ (9a)

2. $k(A\circ B)=(kA)\circ B$, where k is a scalar. (9b)

3. $(A+B)\circ Q=A\circ Q+B\circ Q$ (9c)

4. $A\circ B=E_N^T(A\otimes B)E_M$, where matrix $E_N$ (or $E_M$) is defined as $E_N=[\ e_1\otimes e_1\vdots\cdots\vdots e_N\otimes e_N]$, $e_i=[0\cdots 0\ \underset{i}{1}\ 0\cdots 0]$, $i=1,\cdots,N$, $E_N^T$ is the transpose matrix of $E_N$. $\otimes$ denotes the Kronecker product of matrices. (9d)

5. If A and B are non-negative, then $\lambda_{min}(A)\min\{b_{ii}\}\leq \lambda_j(A\circ B)\leq \lambda_{max}(A)\max\{b_{ii}\}$, where $\lambda$ is the eigenvalue. (9e)

6. $(detA)(detB)\leq det(A\circ B)$, where det( ) denotes the determinant. (9f)

For more details about the Hadamard product see [7].

The weighting residuals of operators p( ), r() and L( ) in equations (8) can be expressed as

$$\int_\Omega p\{\hat{u}\}W_j d\Omega = Ax, \qquad (10a)$$

$$\int_\Omega r\{\hat{u}\}W_j d\Omega = Bx,  \tag{10b}$$

$$\int_\Omega L\{\hat{u}\}W_j d\Omega = Dx,  \tag{10c}$$

where $x = \{c_k\}$, $c_k$'s are the undetermined parameters in equation (3). Therefore, it is intuitively found that the nonlinear formulations of weighted residuals equation (8) can be expressed as

$$(Ax) \circ (Bx) + Dx = b,  \tag{11}$$

where b is the constant vector. For the traditional scheme of weighted residual (equation (7)), the nonlinear formulation can be only expressed in matrix form as [2, 3]

$$D_{n \times n} x + G_{n \times n^2}(x \times x) = b  \tag{12}$$

by using theorem (2.1), where D is given in equation (10c),

$$G_{n \times n^2} = \int_\Omega [p(\phi) \otimes r(\phi)] W_j d\Omega \in C^{n \times n^2}.  \tag{13}$$

The advantages of the formulation (11) over (12) will be discussed in the later section 4. The preceding inferences implicitly assume that the basis functions in equation (3) satisfy all the boundary conditions. However, this is not necessary in general. Among all methods originated from the method of weighted residuals, the more interesting for engineering applications is the FE and Galerkin methods. The residuals of these two methods are weighted by the basis functions in the approximate equation (3), and the boundary conditions can be found by integrating the governing equations by parts, which leads to the so-called weak formulations. By weighting the residuals of governing equation (1) and boundary conditions (2a, b) via basis functions, we have

$$\int_\Omega [p(\hat{u})r(\hat{u}) + L(\hat{u}) - f]\phi_j d\Omega = \int_{\Gamma_2} (q-\bar{q})\phi_j d\Gamma - \int_{\Gamma_1} (u-\bar{u})\frac{\partial \phi_j}{\partial n} d\Gamma \qquad (14)$$

in the standard way. In contrast, we have

$$\int_\Omega p(\hat{u})\phi_j d\Omega \int_\Omega r(\hat{u})\phi_j d\Omega + \int_\Omega [L(\hat{u}) - f]\phi_j d\Omega = \int_{\Gamma_2} (q-\bar{q})\phi_j d\Gamma - \int_{\Gamma_1} (u-\bar{u})\frac{\partial \phi_j}{\partial n} d\Gamma \quad (15)$$

in the present way.

In order to illustrate how the continuity requirements for u can be lowered in weak formulation, consider the heat conduction in a slab with a temperature dependent thermal conductivity [3, 8]

$$\frac{d}{dx}\left[(1 + \alpha T)\frac{dT}{dx}\right] = 0 \qquad (16)$$

or simply

$$\frac{d^2 T}{dx^2} + \alpha T \frac{d^2 T}{dx^2} + \left(\frac{dT}{dx}\right)^2 = 0, \qquad (17)$$

with the boundary conditions

$$T(0) = T_0 \qquad (18a)$$

$$q(L) = q_L, \qquad (18b)$$

where $\alpha$ is a constant. The unknown T is first approximated by $\hat{T}$, namely,

$$\hat{T} = \sum_{k=1}^{n} N_k(x) T_k, \qquad (19)$$

where $N_k(x)$'s are the interpolation basis functions, and $T_k$'s are the unknown nodal temperatures for an one-dimensional line element of finite element method. For the sake of brevity, we only consider the presented innovative scheme of weighting residual here, namely,

$$\alpha \int_{x_1}^{x_n} \left(\frac{d\hat{T}}{dx}\right) N_j dx \int_{x_1}^{x_n} \left(\frac{d\hat{T}}{dx}\right) N_j dx + \alpha \int_{x_1}^{x_n} \hat{T} N_j dx \int_{x_1}^{x_n} \frac{d^2\hat{T}}{dx^2} N_j dx$$

$$+ \int_{x_1}^{x_n} \frac{d^2\hat{T}}{dx^2} N_j dx = \left(\frac{d\hat{T}}{dx} - \overline{q}\right) N_j \bigg|_{x_1}^{x_n} - \left(\hat{T} - T\right) \frac{dN_j}{dx}\bigg|_{x_1}^{x_n} \quad \text{i=1,2,....,n,} \quad (20)$$

where $x_1$ and $x_n$ are the coordinates of the end nodes of certain line element. Using integration by parts to the left-hand side of Eq. (20) and after some deductions, we have

$$\alpha \left[\int_{x_1}^{x_n} \frac{d\hat{T}}{dx} N_j dx\right]^2 - \alpha \int_{x_1}^{x_n} \hat{T} N_j dx \int_{x_1}^{x_n} \frac{d\hat{T}}{dx}\frac{dN_j}{dx} dx + \int_{x_1}^{x_n} \frac{d^2\hat{T}}{dx^2} N_j dx$$

$$+ \alpha \left[\int_{x_1}^{x_n} \hat{T} N_j dx\right]\left(N_j \frac{d\hat{T}}{dx}\bigg|_{x_1}^{x_n}\right) = -\overline{q}\big|_{x_1}^{x_n} - \left(\hat{T} - T\right) \frac{dN_j}{dx}\bigg|_{x_1}^{x_n} \quad \text{i=1,2,....,n} \quad (21)$$

The above formulation can be restated as

$$\alpha(A\hat{T})^{\circ 2} - \alpha(B\hat{T})\circ(H\hat{T}) + D\hat{T} = Q_n - Q_1, \quad (22)$$

where A, B, H and D are the known coefficient matrices, $Q_1$ and $Q_n$ are the constant vectors at two extremities of certain element. The formulation of equations (17) and (18a, b) discretized in the traditional way can be stated as

$$K_{n\times n}\hat{T} + Q_{m\times m^2}(\hat{T}\otimes\hat{T}) = Q_n - Q_1, \quad (23)$$

which is an expression similar to equation (12). K and Q are the known coefficient matrices.

To clarify our idea more general, consider the following example [2, 3, 9]

$$y'' + \sin(y') + 1 = 0; \quad y(0) = 0, y(1) = 1. \quad (24)$$

In terms of the novel weighting residuals scheme, we have

$$\int_0^1 \hat{y}''\phi_j dx + \sin\int_0^1 \hat{y}'\phi_j dx + 1 = 0. \tag{25}$$

Note that here we use the Galerkin scheme with assumption of basis function satisfying boundary conditions. Integration by parts produces

$$-\int_0^1 \hat{y}'\phi_j' dx + \sin\left(-\int_0^1 \hat{y}\phi_j' dx + \hat{y}\phi_j\Big|_0^1\right) + \hat{y}'\phi\Big|_0^1 + 1 = 0. \tag{26}$$

Furthermore, the above equation can be expressed in the Hadamard matrix function form as

$$Dy + \sin^\circ(Ay + a) = b, \tag{27}$$

where D and A are known coefficient matrices, a and b are constant vectors. Similarly, for one more example [2, 3, 9]

$$y' - e^{-y} = 0; \quad y(0)=0, \quad 0 \leq x \leq 1. \tag{28}$$

Nonlinear formulation in the present way can be expressed as

$$Ay - e^{\circ(-B\bar{y})} = b, \tag{29}$$

where $A = -\int_0^1 \hat{y}\phi_j' dx$, $B = \int_0^1 \hat{y}\phi_j dx$, and $b = -\hat{y}\phi_j\Big|_0^1$.

It is worth stressing that for equations (24) and (28), the standard weighting residuals method is by no means capable of expressing its formulation in explicit matrix form. In contrast, we have matrix-form formulations (27) and (29) by using the novel weighted residuals scheme.

As was mentioned earlier, nearly all of popular numerical techniques can be derived from the method of weighted residuals. Only difference among these numerical methods lies in the use of different weightng and basis functions. From the foregoing deduction, it is noted that nonlinear formulations of matrix-separated Hadamard product form can be obtained no matter what weighting and basis functions we use in the presented methodology of weighted residuals. However, as was pointed out by the present author earlier [5 , 6, 7], even if the standard weighted residuals technique is employed, nonlinear formulations in Hadmard product form can be in general directly gotten in the finite difference, finite volume, and collocation (pseudo-spectral) due to the fact that weighting function used in these methods are Dirac function. In other words, integration is actually all circumvented in these numerical techniques. In contrast, only by using the present weighting residuals methodology, the Galerkin, spectral, finite element, boundary element, moments, and least square methods can express their nonlinear formulations in the Hadamard product form.

On the other hand, it is well known that the least square and boundary elements encounter great difficulties to handle nonlinear problems in general. The present innovative scheme of weighting residuals can effectively overcome hindrances in nonlinear computations of these two methods.

For the least square method, it is an uneasy task in the standard way to obtain weight functions of nonlinear operators compared to the linear case. Alternatively, in the present way, we choose the least square weighting functions of each linear operator to separately weight the corresponding linear operator. Consider equation (6), we have

$$\int_\Omega p(\hat{u})\frac{\partial p(\hat{u})}{\partial \phi_j}d\Omega \int_\Omega r(\hat{u})\frac{\partial r(\hat{u})}{\partial \phi_j}d\Omega + \int_\Omega [L(\hat{u})-f]\frac{\partial L(\hat{u})}{\partial \phi_j}d\Omega = 0, \quad j=1,2,...,N, \quad (30)$$

Obviously, the above formulation is absolutely easily done as in the linear case.

For the boundary element method, the deterrence arises from how to eliminate volume integral of nonlinear terms. In terms of the present weighting residuals strategy, we have

$$\int_\Omega p(\hat{u})u_p^* d\Omega \int_\Omega r(\hat{u})u_r^* d\Omega + \int_\Omega [L(\hat{u})-f]u_L^* d\Omega = \int_{\Gamma_2}(q-\bar{q})\phi_j d\Gamma - \int_{\Gamma_1}(u-\bar{u})\frac{\partial \phi_j}{\partial n}d\Gamma, \quad (31)$$

where u*'s with different subscript represent the fundamental solutions relative to each different linear operator. Therefore, the domain integral in lift side of equation (31) can be eliminated easily as in the linear case. However, if p( ), y( ) or even L( ) has not a fundamental solution, it may be feasible to use the recently developed Multiple reciprocity boundary element technique [10] to transfer the domain integral of each separated linear operators in equation (31) to boundary integral in the same way of linear case. This is possibly only a way to apply the multiple reciprocity BEM for general nonlinear problems. We will discuss this aspect further in a separate paper. Also, it may be useful to apply the dual reciprocity BEM to equation (31).

The very basic idea of the weighted residual method is that the weighted average of the error is required to vanish over the solution domain via certain error distribution principle. The philosophy behind the presented methodology of nonlinear MWR is to consider nonlinear operators as a combination of linear operators. Therefore, we apply certain standard MWR technique such as the least square and Galerkin criterions to each of linear operators at first, and then couple these separatedly-weighted linear ones to yield the numerical analogue of corresponding nonlinear operators, finally, set the sum of weighted averages of all nonlinear, linear and constant terms to zero. It is stressed that according to the foregoing discussions, the weighting residual techniques are here understood in a broader way.

In linear computations, some weighting functions are often preferred in some problems (linear operators) or not in other cases. This is because each linear operator no matter which is a linear operator alone or a component of nonlinear operator may has its own distinct error behavior. It is therefore reasonable to employ the different weighting functions which can best reflect error distribution corresponding to each linear operator. We suggest that whenever possible, different weighting functions may be used to handle different linear operator components of complex nonlinear operators. In other words, different numerical technique such as FE, moments and least square methods, etc, can be naturally mixed to deal with one problem simultaneously. For example, weighted residual equation (15) can be restated as

$$\int_{\Omega} p(\hat{u}) W_j^p d\Omega \int_{\Omega} r(\hat{u}) W_j^r d\Omega + \int_{\Omega} [L(\hat{u}) - f] W_j^L d\Omega = \int_{\Gamma_2} (q - \bar{q}) \phi_j d\Gamma - \int_{\Gamma_1} (u - \bar{u}) \frac{\partial \phi_j}{\partial n} d\Gamma \qquad (32)$$

where each different linear operator may use the same or different weight functions. Note that weighting functions in equation (32) is different from those in equation (31) in that they are unnecessarily the fundamental solutions of the corresponding linear operators. In contrast, the standard MWR technique usually employs the same weighting function except in boundary element method, in which weighting functions for governing and boundary equations are different.

The given innovative scheme has the important features of requiring the formulation of all linear differential operators only once, which save considerable computing resources by avoiding the repeated integration of Jacobian matrix in the iterative solution of nonlinear systems. In addition, due to the Hadamard form formulation available in the present WRM nonlinear discretization, simple iteration method (Picard method) can be effectively used to solve these nonlinear algebraic equations as pointed out in Chen [3]. The present strategy is also much simpler to use. Other possible benefits will be given in detail in section 3 including the rapid evaluation of Jacobian matrix via SJT product and nonlinear uncoupling computations.

### 3. SJT product and nonlinear uncoupling computation

The SJT product, introduced by Chen et al. [2, 3, 4], is presented in this section to efficiently calculate accurate Jacobian matrix of nonlinear formulation with the Hadamard product form. A nonlinear decoupling approach given in Chen et al. [3, 4] is extended to nonlinear computation of the presented weighted residuals strategy. These works will show the great advantages of the present novel MRW in computing efficiency.

## 3.1. Caculation of Jacobian matrix by SJT product

The finite difference, collocation methods and their variants such as differential quadrature [11] and pseudo-spectral methods belong to the point-wise approximating numerical technique, which use Dirac weight function in MWR and can express their nonlinear formulation in Hadamard product form [2-4]. In addition, the finite volume, dual reciprocity BEM (the presently only one efficient technique in applying BEM to nonlinear problems) and numerical techniques based on radial basis functions can express their analogue of nonlinear differential operators in the Hadamard product form. However, in their standar form, the FE, Galerkin, moments, BE and least square methods can not have Hadamard product expression of nonlinear analogues. The present novel methodology of weighted residual makes it reality that all of these numerical techniques can formulate their approximate algebraic equations in the Hadamard product form. In the following, the numerical analogues of some nonlinear operators often encountered in practice are given by in the Hadamard product form

- $c(x,y)u_{,x} = \{c(x_j, y_j)\} \circ (A_x U)$,  (33a)

- $(u_{,x})^q = (A_x U)^{\circ q}$, where q is a real number,  (33b)

- $\dfrac{\partial u^m}{\partial x} \dfrac{\partial u^n}{\partial y} = (A_x U^{\circ m}) \circ (A_y U^{\circ n})$,  (33c)

- $\sin u_{,x} = \sin^{\circ}(A_x U)$,  (33d)

- $\exp(u_{,x}) = \exp^{\circ}(A_x U)$.  (33e)

In the above equations $(\ )_{,x} = \partial(\ )/\partial x$; $A_x$ and $A_y$ denote the known coefficient matrices resulting from certain numerical method. From now on we define the nonlinear discretization expression in the Hadamard product, power and function form as the formulation-H. In what follows, the SJT product is introduced to efficiently compute analytical solution of the Jacobian derivative matrix.

**Definition 3.1.** If matrix $A=[a_{ij}]\in C^{N\times M}$, vector $U=\{u_j\}\in C^{N\times 1}$, then $A \lozenge U=[a_{ij}u_i]\in C^{N\times M}$ is defined as the postmultiplying SJT product of matrix A and vector U, where $\lozenge$ represents the SJT product. If M=1, $A \lozenge B = A \circ B$.

**Definition 3.2.** If matrix $A=[a_{ij}]\in C^{N\times M}$, vector $V=\{v_j\}\in C^{M\times 1}$, then $V^T \lozenge A=[a_{ij}v_j]\in C^{N\times M}$ is defined as the SJT premultiplying product of matrix A and vector V.

Considering a simple special case of the nonlinear formulation (33c), we have

$$\frac{\partial}{\partial U}\{(A_xU)\circ(A_yU)\} = A_x \lozenge (A_yU) + A_y \lozenge (A_xU). \tag{34}$$

Formula (34) produces the accurate Jacobian matrix through simple algebraic computations. The SJT premultiplying product is related to the Jacobian matrix of the Hadamard power nonlinear formulations such as $\frac{dU^m}{dx} = AU^{\circ m}$, i.e.,

$$\frac{\partial}{\partial \bar{U}}\{A_xU^{\circ m}\} = \left(mU^{\circ(m-1)}\right)^T \lozenge A_x. \tag{35}$$

In the following, we discuss some operation rules of applying the SJT product to evaluate the Jacobian matrix of the nonlinear formulations (33).

- $\frac{\partial}{\partial \bar{U}}\{\{c(x_j,y_j)\}\circ(A_xU)\} = A_x \lozenge \{c(x_j,y_j)\}$ (36a)

- $$\frac{\partial}{\partial \bar{U}}\left\{(A_x U)^{\circ q}\right\} = q A_x \Diamond (A_x U)^{\circ(q-1)}.$$  (36b)

- $$\frac{\partial}{\partial \bar{U}}\left\{(A_x U^{\circ m}) \circ (A_y U^{\circ n})\right\} = m\left\{(U^{\circ(m-1)}) \Diamond A_x\right\} \Diamond (A_y U^{\circ n}) + n\left\{(U^{\circ(n-1)}) \Diamond A_y\right\} \Diamond (A_x U^{\circ m})$$  (36c)

- $$\frac{\partial}{\partial \bar{U}}\left\{\sin(A_x U)\right\} = A_x \Diamond \cos^{\circ}(A_x U)$$  (36d)

- $$\frac{\partial}{\partial \bar{U}}\left\{\exp^{\circ}(A_x U)\right\} = A_x \Diamond \exp^{\circ}(A_x U)$$  (36e)

- If $\psi = f^{\circ}(\phi)$, $\phi = \varphi^{\circ}(U)$, we have $\frac{\partial \psi}{\partial U} = \frac{\partial \psi}{\partial \phi}\frac{\partial \phi}{\partial U}$.  (36f)

In the above equations $\frac{\partial}{\partial \phi}$ and $\frac{\partial}{\partial U}$ represent the Jacobian derivative matrix of certain Hadamard vector function with respect to vectors $\phi$ and U, respectively. It is observed from these formulas that the Jacobian matrix of the nonlinear formulation-H can be calculated by using the SJT product in the chain rules similar to those in differentiation of a scalar function. The above computing formulas yield the analytical solutions of the Jacobian matrix. The computational effort of a SJT product is only $n^2$ scalar multiplications.

The finite difference method is often employed to calculate the approximate solution of the Jacobian matrix and also requires $O(n^2)$ scalar multiplications. In fact, the SJT product approach requires $n^2$ and $5n^2$ less multiplication operations than one and two order finite differences, respectively. Moreover, the SJT product produces the analytic solution of the Jacobian matrix. In contrast, the approximate Jacobian matrix yielded by the finite difference method affects the accuracy and convergence rate of the Newton-Raphson method, especially for highly nonlinear problems. The efficiency and accuracy of the SJT product approach were numerically demonstrated in [2-4].

### 3.2. Uncoupling Computations

It is usually a difficult task to solve complex coupled nonlinear partial differential equations of large system [12]. Chen and Zhong [3, 4] proposed an uncoupling computation methodology by using the Hadamard product and SJT product. In this section we will introduce concept of the relative Jacobian derivative matrix among different dependent variables to simplify decoupling computations, and discussions are also given on how to apply the Hadamard product and SJT product to uncouple the nonlinear coupling formulation-H's. The following example will help to illustrate and clarify the present technique.

The example is geometrically nonlinear bending of the prismatic beam under uniformly distributed loading. The normalized governing equations are given by[3]

$$\frac{d^4W}{dX^4} - \left(\frac{1}{rl}\frac{dU}{dX} + \frac{1}{2}\left(\frac{dW}{dX}\right)^2\right)\frac{d^2W}{dX^2} = \frac{qL^4}{EIr}, \tag{37}$$

$$\frac{1}{rL}\frac{d^2U}{dX^2} + \frac{dW}{dX}\frac{d^2W}{dX^2} = 0, \tag{38}$$

where A=area of beam cross section, E=modulus of elasticity, I=centroidal moment of inertia of beam cross section, L=length of beam. The variables have been normalized in the form:

$$X = \frac{x}{L}, \quad r^2 = \frac{I}{A}, \quad W = \frac{w}{r}, \quad U = \frac{u}{r}. \tag{39}$$

The corresponding boundary conditions are

$$W = 0, \quad \frac{dW}{dX} = 0, \quad U = 0, \quad \text{at } X=0, \tag{40a}$$

$$M = EI\frac{d^2W}{dX^2} = \overline{M} = 0, \quad Q = -\frac{d^3W}{dX^3} = \overline{Q} = 0, \quad \frac{dU}{dX} = \overline{U} = 0, \, X=1. \tag{40b}$$

In terms of the present MWR method, we have

$$\int_0^1 \left[\frac{d^4W}{dX^4} - \frac{qL^4}{EIr}\right]\phi_j dX - \left\{\frac{1}{rl}\int_0^1 \frac{d^2U}{dX^2}\varphi_j dX + \frac{1}{2}\left[\int_0^1 \frac{dW}{dX}\phi_j dX\right]^2\right\}\int_0^1 \frac{d^2W}{dX^2}\phi_j dX$$
$$= \left[\left(EI\frac{d^3W}{dX^3} + \overline{Q}\right)\phi_j - \left(EI\frac{d^2W}{dX^2} - \overline{M}\right)\frac{d\phi_j}{dx}\right]_{x=1} \tag{41}$$

$$\frac{1}{rl}\int_0^1 \frac{d^2U}{dX^2}\varphi_j dX + \int_0^1 \frac{dW}{dX}\phi_j dX \int_0^1 \frac{d^2W}{d^2X}\phi_j dX = 0, \tag{42}$$

where $\phi$ and $\varphi$ are respectively basis functions of W and U and satisfy the essential boundary conditions. Note that for the brevity, W and U in equations (41) and (42) as well as those afterwards represent the approximate solutions. Furthermore, using integration by parts, we have

$$\int_0^1 \frac{d^2W}{dX^2}\frac{d^2\phi_j}{dX^2}\phi_j dX - \left\{\frac{1}{rl}\left[\frac{dU}{dX}\varphi_j\Big|_0^1 - \int_0^1 \frac{dU}{dX}\frac{d\varphi_j}{dX}dX\right] + \frac{1}{2}\left[\int_0^1 \frac{dW}{dX}\phi_j dX\right]^2\right\}$$
$$\left[\frac{dW}{dX}\phi_j\Big|_0^1 - \int_0^1 \frac{dW}{dX}\frac{d\phi_j}{dX}dX\right] = \int_0^1 \frac{qL^4}{EIr}\phi_j dX + \left[\overline{Q}\phi_j + \overline{M}\frac{d\phi_j}{dx}\right]_{x=1} \tag{43}$$

$$-\frac{1}{rl}\int_0^1 \frac{dU}{dX}\frac{d\varphi_j}{dX}dX + \int_0^1 \frac{dW}{dX}\phi_j dX\left[\frac{dW}{dX}\phi_j\Big|_0^1 - \int_0^1 \frac{dW}{dX}\frac{d\phi_j}{dX}dX\right] = -\frac{1}{rl}\frac{dU}{dX}\varphi_j\Big|_0^1. \tag{44}$$

Finally, we have the formulation-H of this coupled equations

$$DW + \left[(AU + a) + \frac{1}{2}(BW)^{\circ 2}\right]\circ(CW + c) = b, \tag{45}$$

$$AU + (BW)(CW + c) = -a, \tag{46}$$

where D, B and C denote the known coefficient matrices for the desired W, and A is the coefficient matrix for the unknown U. a and b are constant vectors. The unknown vector W is chosen as the basic iterative variable in the present computation. From Eq. (46), we have

$$U = -A^{-1}\left[a + (BW) \circ (CW + c)\right], \tag{47}$$

$$\frac{\partial U}{\partial W} = -A^{-1}\left[B \Diamond (CW + c) + C \Diamond (BW)\right], \tag{48}$$

where

$$\frac{\partial U}{\partial W} = \begin{bmatrix} \frac{\partial U_1}{\partial W_1} & \frac{\partial U_1}{\partial W_2} & \cdots & \frac{\partial U_1}{\partial W_n} \\ \frac{\partial U_2}{\partial W_1} & \frac{\partial U_2}{\partial W_2} & \cdots & \frac{\partial U_2}{\partial W_n} \\ \vdots & \vdots & \ddots & \vdots \\ \frac{\partial U_n}{\partial W_1} & \frac{\partial U_n}{\partial W_2} & \cdots & \frac{\partial U_n}{\partial W_n} \end{bmatrix} \tag{49}$$

is the Jacobian derivative matrix of vector U with respect to vector W. Therefore, the unknown vector U and its Jacobian derivative matrix can be obtained by vector W. The introduction of the relative Jacobian derivative matrix $\frac{\partial U}{\partial W}$ is a key step to simplify uncoupling computations, especially for equations with many coupled dependent variables [3, 4]. We choose equation (45) as the basic equation, namely,

$$\psi(W) = DW + \left[(AU + a) + \frac{1}{2}(BW)^{\circ 2}\right] \circ (CW + c) - b = 0. \tag{50}$$

The corresponding Jacobian matrix can be evaluated by the SJT product, namely,

$$\frac{\partial \psi(W)}{\partial W} = D + \left[A\frac{\partial U}{\partial W} + B \Diamond (BW)\right] \Diamond (CW + c) + C \Diamond \left[(AU + a) + \frac{1}{2}(BW)^{\circ 2}\right]. \tag{51}$$

The resulting Newton-Raphson iteration equation for this case is

$$W^{(k+1)} = W^{(k)} - \left[\frac{\partial \psi(W^k)}{\partial W}\right]^{-1} \psi(W^{(k)}). \tag{52}$$

It is seen from the above manipulations that the coupled partial differential equations of two variables are decoupled into the single variable equation. Therefore, the computational effort and storage requirements are reduced to about only one-eighth and one-fourth, respectively, similar to what we have done for complex coupled von Karman equations of geometrically nonlinear bending plates [4]. Also, the repeated integration of Jacobian matrix is no longer needed in the present procedure.

## 4. Remarks

One of the major factors which affects the efficiency of the FE and Calerkin methods is need to repeat numerical integration of Jacobian (stiffness) matrices. The presented formulation cures this deficiency and yet they maintain strong geometric and boundary flexibility with the weak formulation. The conventional nonlinear FE method may also be too complex mathematically for routine applications. In contrast, the present Hadamard product formulation is an explicit and simple matrix approach. By using the novel weighted residual formulation, the least square and boundary element methods can be conveniently applied to nonlinear analysis.

The universal formulation-H provides a computational attractiveness to develop the unified techniques in the nonlinear analysis and computation. The SJT product and the relative uncoupling algorithms given in section 3 are the first successful attempts in this aspect. References [2-4] provided some benchmark numerical examples to demonstrate the superiority of these approaches. More numerical experiments of the present strategy will be provided in a sequent paper.